\documentclass[a4paper,12pt]{amsart}
\usepackage{tikz}
\usepackage{tikz-cd}
\usetikzlibrary{matrix, arrows}
\usepackage{tikz}

\usepackage{bm}
\usepackage{amsfonts,latexsym,rawfonts,amsmath,amssymb,amsthm, mathrsfs, lscape,color}
\usepackage{extsizes}
\usepackage[all]{xy}
\usepackage[english]{babel}		\usepackage[T1]{fontenc}
\usepackage{graphicx}
\usepackage{booktabs}
\usepackage{array, tabularx}
\usepackage{setspace}
\usepackage{textcomp}

\usepackage{tikz}
\usepackage{multirow}
\usepackage{paralist}
\usepackage{comment}

\usepackage[hypertexnames=false,
backref=page,
    pdftex,
    pdfpagemode=UseNone,
    breaklinks=true,
    extension=pdf,
    colorlinks=true,
    linkcolor=blue,
    citecolor=blue,
    urlcolor=blue,
]{hyperref}

\usepackage[top=1.5in, bottom=.9in, left=0.9in, right=0.9in]{geometry}

\newcommand{\referenza}{}

\newtheorem{thm}{Theorem}[section]
\newtheorem*{thm*}{Theorem \referenza}
\newtheorem{cor}[thm]{Corollary}
\newtheorem*{cor*}{Corollary \referenza}
\newtheorem{lem}[thm]{Lemma}
\newtheorem*{lem*}{Lemma \referenza}

\newtheorem{prop}[thm]{Proposition}
\newtheorem*{prop*}{Proposition \referenza}
\newtheorem{prob}[thm]{Problem}

\newtheorem*{conj*}{Conjecture \referenza}
\newtheorem{rmk}[thm]{Remark}
\newtheorem*{rmk*}{Remark}

\newtheorem{exa}[thm]{Example}
\newtheorem{defi}[thm]{Definition}
\newtheorem{exe}[thm]{Exercise}

\def\bexc{\begin{exe}}
\def\eexc{\end{exe}}

\def\be{\begin{equation}}
\def\ee{\end{equation}}
\def\bp*{\begin{prop*}}
\def\ep*{\end{prop*}}
\def\bt*{\begin{thm*}}
\def\et*{\end{thm*}}
\def\bt{\begin{thm}}
\def\et{\end{thm}}
\def\bp*{\begin{prop*}}
\def\ep*{\end{prop*}}
\def\bp{\begin{prop}}
\def\ep{\end{prop}}
\def\bl*{\begin{lem*}}
\def\el*{\end{lem*}}
\def\bl{\begin{lem}}
\def\el{\end{lem}}
\def\bec*{\begin{cor*}}
\def\enc*{\end{cor*}}
\def\bc{\begin{cor}}
\def\ec{\end{cor}}
\def\bex{\begin{exa}}
\def\eex{\end{exa}}
\def\de{\begin{defi}}
\def\ed{\end{defi}}
\def\br{\begin{rmk}}
\def\er{\end{rmk}}
\def\pb{\begin{prob}}
\def\eb{\end{prob}}

\def\k{\Bbbbk}

	\def\tens{\otimes}

\newcommand{\lra}{\longrightarrow}

\def\hra{\hookrightarrow}
  
\def\ES{\varnothing}
\numberwithin{equation}{section}

 \def \t\G {\widetilde \Gamma} \def \t\T {\widetilde \T}
\def\bit{\begin{itemize}} \def\eit{\end{itemize}}
\def\benu{\begin{enumerate}} \def\enu{\end{enumerate}}
\def\demo{\begin{proof}} \def\enddemo{\end{proof}}

\def \p {\partial}
\def\de{\p}

\numberwithin{equation}{section}

\def\acal{\mathcal A}  \def\ccal{\mathcal C} \def\dcal{\mathcal D}\def\ecal{\mathcal E} \def\fcal{\mathcal F}
\def\hcal{\mathcal H}\def\ical{\mathcal I}
\def\ncal{\mathcal N}\def\ocal{\mathcal O}\def\rcal{\mathcal R}
\def\tcal{\mathcal T}
    
\def\imsc{\mathscr I}

     \def\Smfr{\mathfrak{S}}

   \def\a {\alpha} \def\b {\beta}
\def\N{{\mathbb N}}     \def\d {\delta} 
\def\C{{\mathbb C}}      \def\l{\lambda}
					\def\p{\partial}
\def\r{\varrho} 
					\def\s{\sigma}\def\t{\theta}\def\z{\zeta}
\def\P{{\mathbb P}}

\def\T{{\mathbb T}}

\def\V{{\mathbb V}}

 \def\G{\Gamma}

   \def\a {\alpha} \def\b {\beta}
 
\def\L{\Lambda}
\def\N{{\mathbb N}}     \def\d {\delta}

\def\C{{\mathbb C}}

\def\G{{\mathbb G}}

    \def\l{\lambda}

					\def\p{\partial}
\def\r{\varrho}
					\def\s{\sigma}\def\t{\theta}\def\z{\zeta}
\def\P{{\mathbb P}}
\def\T{{\mathbb T}}

\def\V{{\mathbb V}}

   \def\a {\alpha} \def\b {\beta}
\def\N{{\mathbb N}}     \def\d {\delta} 
\def\C{{\mathbb C}}      \def\l{\lambda}
				\def\p{\partial}
\def\r{\varrho}
				\def\s{\sigma}\def\t{\theta}\def\z{\zeta}
\def\fa{\forall}

 \def\G{\Gamma}

\def\1{\rm 1} \def\2{\rm 2} \def\3{\rm 3} \def\4{\rm 4} \def\5{\rm 5} \def\6{\rm 6}

\def\lp{\rm(\,} \def\rp{\,\rm)}
\def\smi{\setminus} \def\ssmi{\!\smallsetminus\!}

\def\k{\Bbbk}
\def\Til {\widetilde}
\def \Hat {\widehat}

\def\IN{\infty}

\def\kdim{{{\dim}_\Bbbk}\,}

\def\tms{\times}

 \def\oli{\overline}
\def\sbs{\subset} \def\sps{\supset}
\def\sbseq{\subseteq}  
\def\nsbs{\nsubseteq}

\def\beqn{\begin{eqnarray*}} 
\def\eeqn{\end{eqnarray*}}
\def\beqnn{\begin{eqnarray}} 
\def\eeqnn{\end{eqnarray}}
\def\ba{\begin{aligned}}
\def\ea{\end{aligned}}
\def\bca{\begin{cases}}
\def\eca{\end{cases}}
\def\lra{\longrightarrow} 
\def\neqv{\not\equiv}
\def\rsa{\rightsquigarrow} 
\def\nin{\noindent}
\def\bit{\begin{itemize}}
\def\eit{\end{itemize}}

\allowdisplaybreaks[1]

\title[Convexity and Concavity from Complex Variables to Algebraic Geometry] 
{Convexity and Concavity from Complex Variables to Algebraic Geometry}

\author{ Giuseppe Tomassini}
\address[]{Scuola Normale Superiore, Piazza dei Cavalieri, 7 - I-56126 Pisa, Italy
}
\email{giuseppe.tomassini@sns.it}
   
\thanks{}

\date{\today}

\begin{document}
\baselineskip=.6 cm


\maketitle
\tableofcontents
\section{Introduction}\label{gen}
\nin The topics {\em Convexity and concavity}\footnote{Observe that the first  notion of (algebraic) convexity goes back to Goodman and Landman\cite{GL}}
and {\em Envelopes} are central in Complex analysis and extensively investigated. The aim of this paper  is to find a possible counterpart in Algebraic geometry. The article presents preliminary results on this topic, laid out in five sections. In Section 2 we discuss some elementary results on the extensions of morphisms, while in Section 3 we consider affineness and envelopes. The main result is Theorem \ref{p8} on the existence of basis of affine neighborhoods of an affine subvariety. In the same section,  we introduce the notion of {\em affine envelope} of an algebraic variety and the main steps of a path to prove the existence of the affine envelope of a non-singular algebraic variety $X$ endowed with an étale morphism $X\to D$ onto an open subset $D$ of $\k^n$, $n\ge 2$, Theorem \ref{env1}. In Section 4, we provide the notions of {\em strongly $\ocal$-convexity} and {\em strongly $\ocal$-concavity} of an algebraic variety. The main result of the section is Theorem \ref{p35} which compares semi-affineness and strongly $\ocal$-convexity. The section also contains a formulation of the {\em algebraic Levi problem}. Finally, in the last section, we give the notions of  {\em convexity and concavity along the fibers} of a line bundle. As an application of Theorem \ref{p35}, we prove two existence theorems of modifications (Theorem \ref{t6} and Theorem \ref{t4}).

It is worth observing that the considerations made so far can be carried out for more general classes of {\em algebraic objects} (e.g. linear fiber bundles $\sf V(\fcal)$, cohomology groups $H^r(X,\fcal)$, where $\fcal$ is a coherent sheaf). We address these cases in a work in preparation.

I would like to thank Giorgio Ottaviani for useful discussions I had with him.
\section{Some elementary results}
\subsection{Notations}\label{gen}For a given algebraic variety\footnote{In what follows, by algebraic variety we intend a $\k$-algebraic variety in the sense of Serre where $k$ is an algebraically closed field of characteristic 0.} $X$ we denote ${\ocal}(X)$ the $\Bbbk$-algebra of the regular functions and ${\mathcal R}(X)$ the function field of $X$. 
 $X$ is said {\em semi-affine} if it is proper over an affine variety $Y$ (which is then called an {\em affine reduction} of $X$) and {\em quasi-affine} if it is an open subset of an affine variety. 

Let $\mathfrak S(X)$ be the set of the maximal ideals of $\ocal(X)$. If $X$ is affine, $\mathfrak S(X)$ is nothing but the set $\{M_x\}_{x\in X}$ of all maximal ideals associated to points of $X$. In the general case for a maximal ideal $M\in\Smfr(X)$ one has
$$
\k\sbs\ocal(X)/M\sbs \k(x_1\ldots,x_n).
$$
If $\ocal(X)/M$ is an algebraic extension of $\k$, then $\k=\ocal(X)/M$. 

We denote ${\mathfrak S}_0(X)\sbs \mathfrak S(X)$ the set of those maximal ideals $M$ such that $\ocal(X)/M\simeq \k$. $\mathfrak S_0(X)(X)$ contains the set $\{M_x\}_{x\in X}$.

The algebra $\ocal(X)$ determines a function algebra $\Hat{\ocal(X)}$ on $\Smfr_0(X)$ in the following way: for every $f\in\ocal(X)$, $\Hat{f}$ is the function which associates to $M\in\Smfr_0(X)$ the class of $f$ in $\ocal(X)/M$. Clearly, the map $\ocal(X)\to\Bbbk$, $f\mapsto{\Hat f}(M)$, is a character $\chi$ of $\ocal(X)$ and $M=\chi^{-1}(0)$. Conversely, every character $\chi$ of $\ocal(X)$ is determined by the maximal ideal $M=h^{-1}(0)$. In particular, if $x\in X$ and $f\in\ocal(X)$ then $f(x)=\Hat f(M_x).$ In what follows we identify liberally characters $\chi$ of $\ocal(X)$ and maximal ideals $h^{-1}(0)\in\Smfr_0(X)$. Accordingly, for every $M=\chi^{-1}(0)$ we set $\Hat f(\chi)=\Hat f(M)$. 

We consider on $\Smfr_0(X)$ the weak topology induced by $\Hat{\ocal(X)}$. Let $M_0\in\Smfr_0(X)$, $\Hat f_1,\ldots,\Hat f_n$ in $\Hat{\ocal(X)}$ and a finite subset $F\sbs\k$ such that $\Hat f_1(M_0),\ldots,\Hat f_n(M_0)\notin F$. As $\Hat f_1 \ldots,\Hat f_n$ vary in $\Hat{\ocal(X)}$ and $F$ in $\k$ the subsets
$$
\{M\in\Smfr_0(X):\Hat f_j(M)\notin F, 1\le j\le n\}
$$
form a basis of neighborhoods of $M_0$ for the weak topology.

The natural map $\d:X\to\Smfr_0(X)$, $x\mapsto M_x$ is continuous and not surjective in general (e.g. $X=\k^2\ssmi\{{\sf 0}\}$).
\subsection{Weakly regular functions} Let $A$ be a (closed) subvariety of an algebraic variety $X$ and $\mathscr{I}_A$ the ideal of $A$. 
 A regular function $f:X\ssmi\, A\to\k$ is said {\em $\mathscr{I}_A$-dominated {\lp}along $A${\rp}} if, for every point $z\in A$, $h\in\mathscr{I}_{A,z}$ and a neighborhood $V$ of $0\in\k$ there exists a neighborhood $U$ of $z$ in $X$ such that $hf(U\ssmi A)\sbs V$. In other words if
\be\label{Rie}
\lim\limits_{\stackrel{x\to A}{x\in X\ssmi A}}h(x)f(x)=0.
\ee
A $\mathscr{I}_A$-dominated function $f:X\ssmi\, A\to\k$ is called {\em weakly regular} in $X$. 
We have the following
\bl
Assume that $X$ is normal. Then every weakly regular function $f:X\ssmi A\to\k$ extends regularly to $X$. 
\el
\demo
Since $X$ is normal it is enough to show the lemma when $X$ and $A$ are both non-singular and $A$ is of codimension 1. Let $\imsc=\imsc_A$ be the ideal of $A$ and take for $h$ the generator of $\imsc_a$, $a\in A$. Then $h$ generates $\imsc_x$ for $x$ near $a$. Let $f\in \ocal(X\ssmi A)$ satisfying \ref{Rie}: $f$ is a rational function in $X$, $f=P/Q$ with $P$, $Q$ regular near $a$. If $Q\neqv 0$ on $U\cap A$ then, since $X$ is regular and $\{Q=0\}\cap A$ is of codimension $\ge 2$, $f$ extends through $A$ (if $X$ is a curve then  $Q(a)\neq 0$ and $f$ is regular).

Let $Q\equiv 0$. Then, since $X$ is non-singular in particular locally factorial, near $a$ we have $Q=\l h^m$ with $m\ge 0$ minimum, i.e. $P\neqv 0$ on $A$. Then $f=P/\l h^m$. If $m>0$, the condition \ref{Rie} implies 
$$
\lim\limits_{\stackrel{x\to A}{x\in X\ssmi A}}P(x)=\lim\limits_{\stackrel{x\to A}{x\in X\ssmi A}} h(x)^mf(x)=0
$$
which is absurd since $P\neqv 0$ on $A$. Hence $m=0$ and $f$ is regular. 
\enddemo
\br
Using the normalization, it is easy to show that weakly regular functions are integral over the ring of the regular functions. The converse is also true thanks to the existence of the universal denominator. 
\er
\subsection{Extension of morphisms}\label{gaho1}
 \subsubsection{Gaps and holes} Let $X,Y$ be algebraic varieties, $A\neqv X$ a closed subvariety and
$f:X\ssmi A\to Y$ a morphism. Let 
$$
\Gamma(f)=\big\{(x,f(x))\in X\tms Y:x\in X\ssmi A\big\}
$$ 
be the graph of $f$ and $\oli{\Gamma(f)}$ its closure.

For every subset $S$ of $X$ 
$$
f(S):={\rm pr}_Y\big((S\tms Y)\cap\oli{\Gamma(f)}\big)
$$
is, by definition, the {\em image} of $S$ under $f$. We set $f(x):=f(\{x\})$. 

 A point $x\in X$ is said a {\em hole} of $f$ if $x\notin{\rm pr}_X(\oli{\Gamma(f)})$ and it is said a {\em gap} if either $f(x)=\ES$ or $f(x)\neq\ES$ but $f(x)$ is not complete. 

A morphism $f:X\to Y$ is said {\em gapless} if $f(x)$ is a non-empty, complete subvariety 
of $Y$ for every $x\in X$. As it is easily seen, $f$ is gapless if ${{\rm pr}_X}_{|\oli{\Gamma(f)}}$ is a proper morphism. Trivially, gapless implies no hole.
\bl\label{gaho}
Let $X,Y$ be algebraic varietes, $A\neqv X$ a closed subvariety and $f:X\ssmi A\to Y$ a morphism. If $x\in A$ and $y\in f(x)$ there exists a curve $C\sbs X$ near $x$ such that $C\cap A=\{x\}$ and $y\in f(C)\cap f(A)$. 
\el 
\demo
We may assume that $Y$ is complete.
Consider a desingularization $p:\G^\ast\to\oli{\G(f)}$ of $\oli{\G(f)}$ and a point $w^\ast\in\G^\ast$ such that $p(w^\ast)=(x,y)$. Let $U$ be an open neighborhood of $w*$ and $C^\ast\sbs U$ be a closed curve through $w^\ast$ such that 
$$
C^\ast\cap p^{-1}(\oli{\G(f)}\ssmi\G(f))=w^\ast.
$$
Since ${\rm pr}_X\circ p$ is proper there exists a neighborhood $V$ of $x$ in $X$ such that $U':=({\rm pr}_X\circ p)^{-1}(V)\sbs U$ and ${\rm pr}_X\circ p$ is proper from $U'$ to $V$. Then the projection $({\rm pr}_X\circ p)(C_{Uì}^\ast)$ has the desired properties.
\enddemo
In the sequel we will use the Zariski's Main Theorem (ZMT) in the following form:
\bt
Let $X$ be a normal variety, $A$ a closed subvariety with ${\rm  codim} \,A_x \ge 1$. Let $f$ be a morphism $X\ssmi A\to Y$. Then for a point $x\in A$ one of the following instances occurs:
\bit
\item[i\rp] $f(x)=\ES$;
\item[ii\rp] $f(x)$ is a point and then $f$ is regular at $x$; 
\item[iii\rp] $f(x)$ is connected of positive dimension at every point.
\eit
\et
%
As a consequence, the set ${\rm Sing}(f)$ of the non-regular points of $f:X\ssmi A\to Z$ is the algebraic subset 
$$
\big\{x\in X:\kdim f(x)>0\big\}.
$$
\bp
Under the hypothesis of the above theorem if $f$ has no hole then
$$
\kdim {\rm Sing}(f)_x\le\kdim X-2
$$
for every $x\in X$.
\ep
\demo
We may assume that $Y$ is complete. Let $n=\kdim X$ and 
$$
r(x,y)=n-\kdim f(x)_y
$$
be the rank of the regular map ${\rm pr}_X:\Gamma(f)\to X$.
Since $f$ has no hole 
\beqn
{\rm Sing}(f)&=&\big\{x\in X:\kdim f(x)>0\big\}=\\
&&{\rm pr}_X\big\{(x,y)\in X\tms Y:r(x,y)<n\big\}
\eeqn
$$
S:=\big\{(x,y)\in X\tms Y:r(x,y)<n\big\}
$$ 
is an algebraic set and ${{\rm pr}^{-1}_X(x)}\hspace{-2mm}\sbs \hspace{-2mm}S$ for every $x\in{\rm Sing}(f)$. Moreover, since $\kdim f(x)_y\ge1$, for every $(x,y)\in S$ the rank of ${\rm pr}_X|S$ at a point $(x,y)$ is
$$
r_S(x,y)=\kdim S_{(x,y)}-\kdim   
f(x)_y\le n-2.
$$
This implies  
$$
\kdim {\rm Sing}(f)_x\le\kdim X-2.
$$
\enddemo
\subsubsection{Extension}
An algebraic variety $X$ is said to be $\rm K$-{\em complete} if for a given point $y\in X$ there are finitelely many regular functions $f_1,\ldots,f_m\in\ocal(X)$ such that $y$ is an isolated point
of  $$ \{x\in X:f_j(x)=f_j(y), 1\le j\le m\}. 
$$
An affine variety is $\rm K$-complete. 
\bt\label{p1}
Let $X,Y$ be algebraic varieties with $X$ normal. Let $A\sbs X$ be a subvariety with ${\rm codim} \,A_x \ge 2$, for every $x\in A$ and $f:X\ssmi A\to Y$ a morphism. If $Y$ is affine or $\rm K$-complete and $f$ is without holes, then $f$ extends to a morphism $X\to Y$.
\et 
\demo 
If $Y$ is affine, then $f=(f_1,\ldots,f_n)$, $f_j\in\ocal(X)$, $1\le j\le n$. Since $X$ is normal the thesis is a consequence of the extension theorem for regular functions.
Assume that $Y$ is $\rm K$-complete and that $f$ has no hole. Let $x_0\in A$. Then $\Sigma_f(x_0)\neqv\ES$. Let $y_0\in\Sigma_f(x_0)$ and $h$ be a regular function on $Y$; then, since $X$ is normal, $h\circ f$ extends on all of $X$ by a regular function $g$. By Lemma \ref{gaho} there is a net $\{x_\l\}_{\l\in\Lambda}$ with the properties $x_\l\to x_0$ and $f(x_\l)\to
y_0$. Then $g(x_0)=h(y_0)$ and consequently
$$
\Sigma_f(x_0)\sbseq\bigcap\limits_{h\in\mathcal O(Y)}\big\{y\in Y:h(y)=g(x_0)=h(y_0)\big\}=\{y_0\}
$$
since $Y$ is $\rm K$-complete. Thus for every $x_0\in A$, $\Sigma_f(x_0)$ is a point. Then, ZMT implies that the map $x\mapsto\Sigma_f(x)$ is the desired extension of $f$.
\enddemo
                                                                                           \subsubsection{The complex case} If $\k=\C$ the extension of regular maps of algebraic varieties reduces to the corresponding problem for holomorphic maps between $(X^{\rm an},{\mathcal O}_{X^{\rm an}})$ the complex space associated to a complex algebraic variety $(X,{\mathcal O}_X)$. Then we have the following 
\bp\label{anal1}
Let $X$, $Z$ be complex algebraic varieties, $A\sbs X$ an algebraic subvariety of positive codimension at every point. A holomorphic map $f:X^{\rm an}\to Z^{\rm an}$ which is regular on $X\ssmi A$ is regular. 
\ep
\demo
The problem is local so we may assume that $X\sbs\C^n$ and $Z=\C$.
We first consider the case when $X$ is irreducible. Let $x\in A$, $\imsc_x$ be the ideal of $X_x$ and $I_x$ the ideal of $X^{\rm an}_x$: $I_x=\imsc_x\mathcal O^{\rm an}_{\C^n,x}$ (\cite{Se}). Since $f_{|X\setminus A}$ is regular it extends on $X$ as a rational function, i.e. there exist two polynomials $p,q\in\C[z_1,\ldots,z_n]$, $q\notin \imsc_x$ such that $q f-p=0$ on $X$. If $q(x)\neq 0$ then $f\in \mathcal O_{X,x}$, i.e. $f$ is regular at $x$. Suppose, at contrary, $q(x)=0$. Let $\imsc'_x$ be the ideal generated by $\imsc_x$, $q$ and $I'_x=\imsc _x\mathcal O^{\rm an}_{\C^n,x}$. We have $I'_x=I_x+q\,\mathcal O^{\rm an}_{\C^n,x}$ and $p\in I'_x$. Since $\mathcal O^{\,\rm an}_{\C^n,x}$ is faithfull flat over $\mathcal O_{\C^n,x}$ (see \cite{Se}), $\imsc'_x=I'_x\cap\mathcal O_{\C^n,x}$, hence $p\in \imsc '_x$. It follows that there exists a polynomial $g$ such that $p-gq\in\imsc_x$, therefore $p_{|X}=g_{|X}q_{|X}$ and $f(x)=g(x)$ if $q(x)\neq 0$. Since $X$ is irreducible, $f=g$ everywhere.

In the general case let $X_1\cup X_2\cdots X_m$ be the decomposition of $X$ into irreducible components. We argue by induction on $m$. Set $Y=X_1\cup X_2\cdots X_{m-1}$ and consider the Mayer-Vietoris sequence applied to $Y\cup X_m$
\be\label{an2}
0\to\mathcal O_X\stackrel{\r}\to\mathcal O_Y\oplus\mathcal O_{X_m}\stackrel{\d}\to\mathcal O_{Y\cap {X_m}}\to 0
\ee
where $\r(g)=g_{\vert Y}\oplus g_{|X_m}$ and $\d(g\oplus h)=g_{|Y}-h_{|X_m}$. Since $X$ is affine, the associated cohomology sequence 

\be\label{an3}
\xymatrix{{\ocal_X(X)}\ar[r]^{\hspace{-2cm}\r}&{\ocal _Y(Y)\oplus\ocal_{X_m}(X_m)}\ar[dl]_{\d}\\
\ocal_{Y\cap {X_m}}(Y\cap {X_m})\ar[r]^{}&0}
\ee
is exact.
Let $f$ be like in the statement. Then, by what already proved, $f_{|X_m}$ is regular and by the induction hypothesis $f_{|Y}$ is also regular. Moreover, $f_{|X_m}-f_{|Y}=0$ on $Y\cap{X_m}$. Then, from \ref{an3} it follows that $f$ is regular.
\enddemo
\vskip1.5truecm
\section{Affineness}
\bp\label{pxx}
Let $X$ be an algebraic variety, $Y\sbs X$ a subvariety. If $X\ssmi Y$ is semi-affine, then $Y$ is of pure codimension $1$. Conversely, let $X$ be  affine and locally factorial\footnote{i.e. all the local rings $\ocal_{X,x}$ are UFD.} and $Y$ of pure codimension $1$, then,  $X\ssmi Y$ is affine.
\ep
\demo
We may assume that $X$ is normal. Let $X\ssmi Y$ be semi-affine and $\pi :X\ssmi Y\to U$ a proper morphism where $U$ is an affine variety. Set $Y=Y'\cup Y''$ where $Y'$ is the union of all irreducible components of $Y$ of codimension $1$. Assume, by a contradiction, that $y\in Y''\neqv\ES$. Since $X$ is normal, $\pi$ extends through $Y''$. Let $a=\pi(y)$ and $C\sbs U$ be an algebraic curve through $a$. By construction, in the diagram 
\be
{\begin{split}
\xymatrix{(X\ssmi Y)\tms_U(C\smi\{a\})\ar[r]\ar[d]^{\pi\tms_U {\sf i}}& X\ssmi Y\ar[d]^\pi\\
C\ssmi\{
a\}\ar[r]^{\sf i} & U}
\end{split}}
\ee
the morphism $\p\tms_U {\sf i}$ is not closed and this contradicts the properness of $\pi$.

Conversely let $X$ be locally factorial and affine. Since $Y$ is of pure codimension $1$ and $X$ is locally factorial, $Y$ is locally the zero set of a regular function, then $X\ssmi Y$ is locally affine whence it is affine since $X$ is affine (\cite{Haa}, p. 128).
\enddemo
\br
In the second part of the Proposition \ref{pxx} the hypothesis ``$X$ is affine'' cannot be dropped (e.g. take for $X$ the blow-up of $\Bbbk^2$ at a point). 
\er
\br\label{rem}
Let $X$ be a complex algebraic variety and $(X^{\rm an},\hcal)$ the complex space associated to $X$. Denote $ \hcal(X^{\rm an})$ the algebra of the holomorphic functions in $X^{\rm an}$. If $X$ is affine, then $X^{\rm an}$ is Stein. The converse is not true {\lp\cite[Section 7]{Ne}, see also \cite[Section 3]{Ha}\rp}. Observe that in the counterexample of \cite{Ne} the algebraic variety has no regular global function except the constants. This suggests the following question: is a complex algebraic variety $X$ affine if $\ocal(X)$ separates points and $X^{\rm an}$ is Stein?
 \er
\subsubsection{Affine neighborhoods}
 In complex analysis the following is true: every Stein subspace $Y$ of a complex space $X$ has a basis of Stein neighborhoods (cfr.\cite{Si}; see also \cite{De} for a generalization to $q$-complete spaces). The same question can be raised for algebraic varieties: does a (non necessarily closed) affine subvariety $Y$ of an arbitrary algebraic variety $X$ have a basis of affine neighborhoods? The answer is negative in general. Here is an example.
\bex\label{ex}
Let ${\ccal}\sbs\P^2$ be a regular cubic (\cite{Haa}, \cite{Mi}). Let ${\rm Pic}(\ccal)$ be the group of the Weyl divisors $\sum_{j=1}^m P_j$ of $C$ , ${\rm Pic}_0(\ccal)$
the subgroup of those of degree 0. We have the following theorem: If $P_0$ is a fixed inflection point of $\ccal$ then  $3(P-P_0)=0$ in ${\rm Pic}(\ccal)$ if and only if $P$ is an inflection point. Take a point $P$ such that $P-P_0$ is  no torsion in ${\rm Pic}(\ccal)$. Then no affine open subset of $U\sbs\P^2$ can contain the affine curve ${\ccal}\smi\{P\}$. In fact, if this were the case, there would exist a projective curve ${\dcal}\sbs\P^2$  such that ${\dcal}\cap\ccal=\{P\}$. Then, if  H is the divisor of the hyperplane section, since $P_0$ is an inflection point, we derive $H=3P_0$, $nP=mH$ with $m,n\in\N$ and consequently 
$$
0=nP-mH=3m(P-P_0)
$$
i.e. $P-P_0$ is torsion in ${\rm Pic}(\ccal)$: contradiction.
\eex
\begin{prob}
It could be interesting to eneralize the above example.
\end{prob}
The following theorem provides a positive answer to the question:
\bt\label{p8}
Let $X$ be an algebraic variety $Y$ a (non necessarily closed) affine subvariety $Y$. Then $y$ has a basis of affine neighborhoods if and only if it is closed in an affine open subset of $X$.
\et
\demo
 If $Y$ has a basis of affine neighborhoods then it is closed in an affine open subset of $X$. To prove the converse we may assume that $X$ is affine, $X\sbs \Bbbk^N$, and $Y$ is closed in $X$. Let $U=X\ssmi A$ be a neighborhood of $Y$, where $A$ is the common zeros set of the polynomials $F_1,F_2,\ldots,F_n$. Let $Y_1,Y_2\cdots Y_m$ be the irreducible components of $Y$ and fix points $y_1\in Y_1,y_2\in Y_2,\ldots,y_m\in Y_m$ such that  
$$
y_\nu\in Y_\nu\smi Y_1\cup Y_2\cup\cdots\cup Y_{\nu-1}\cup Y_{\nu+1}\cup\cdots\cup Y_m
$$
for $1\le\nu\le m$.

\nin Then there exists $\l=(\l_1,\ldots,\l_n)\in\Bbbk^n$ such that the polynomial $G_\l:=\sum_{j=1}^n\l_jF_j$ does not vanish at $y_1,y_2,\ldots,y_m$. It follows that $Z=\{G_\l=0\}\cap Y$ is a proper algebraic subset of $Y$ with $\kdim Z_y
<\kdim Y_y$ for every $y\in Z$ (and $A\sbs Z$).
Let $Z_1,Z_2\cdots Z_k$ be the irreducible components of $Z$. Arguing as above we find a polynomial $G_\mu:=\sum_{j=1}^n\mu_jF_j$ such that $W:=\{G_\mu=0\}\cap Z$ is a proper subvariety of $Z$ with $\kdim W_y
<\kdim Z_y$ for every $y\in W$. After many but finite steps we find polynomials $G_1,\ldots,G_p$ with no common zeros on $Y$ and such that 
$$
A\sbs\bigcap_{1\le j\le p}\{G_j=0\}.
$$ 
\nin Since $Y$ is affine and $G_1,\ldots,G_p$ have no common zeros on $Y$, by Theorems A and B there exist polynomials $Q_1,\ldots,Q_p$ such that $\sum_{j=1}^p(Q_jG_j)_{|Y}=1$. It follows that the hypersurface $S=\{\sum_{j=1}^pQ_jG_j)=0\}$ contains $A$ and it is disjoint from $Y$™. It follows that $V=X\ssmi S\sbs U$ is an affine neighborhood of $Y$.      
\enddemo

\subsection{Affine envelopes}\label{afen}
Let $X$ be an algebraic variety. An {\em affine envelope} of $X$ is a pair $(\Til X,{\sf j})$ where $\Til X$ is an affine variety and ${\sf j}$ is an open embedding $X\hra \Hat X$ such that ${\sf j}^*:\ocal(\Til X)\lra\ocal(X)$ is an isomorphism. 
\bp\label{env}			Every open subset $U$ of a locally factorial affine variety $X$ has an affine envelope which is still an open subset of $X$. 
\ep
\demo
Let $U=X\ssmi Y=X\ssmi A\ssmi B$, where $A$ is the union of all irreducible components of $Y$ of pure codimension $1$. Since $X$ is locally factorial, in particular locally normal, every regular function in $U$ extends to $X\ssmi A$. If $A=\ES$ then $X$ is the envelope of $U$. Suppose $A\neqv\ES$. Then, since  $X$ is locally factorial and $A$ is of pure codimension $1$, $X\ssmi\, A$ is locally affine, hence affine (\cite{Haa}, p. 128).
\enddemo
The affine envelope is uniquely determined up to isomorphisms, provided it exists. In such a situation $(\Til X,{\sf j})$, as set of points, is the maximal spectrum ${\Smfr}\,\ocal(X)$ of $\ocal(X)$.
In particular, if $(\Til X,{\sf j})$ exists then $\ocal(X)$ separates points and gives local coordinates.

Let ${\sf V ar}/\k$ be the category of the algebraic varieties, ${\sf F}_X$ the contravariant functor 
\be\label{funct}
{\sf V ar}/\k\lra\{\k-{\sf algebras}\}, 
\ee
defined by
\be\label{fun1}
Y\rsa {\rm Hom_{alg}}(\ocal(X),\ocal(Y)).
\ee
${\sf F}_X$ is representable in ${\V ar}/\k$ if and only if $X$ has an affine envelop (see Section \ref{cxcv} below).

It should be possible to prove the following: 
\bt\label{env1}
Let $X$ be a non-singular algebraic variety, $p:X\to D$ an étale morphism onto an open subset $D$ of $\k^n$, $n\ge 2$. Let ${\sf j}:D\hra\Til D$ be the affine envelope of $D$. Then, the spectrum $\Smfr(X)$ of $\ocal(X)$ can be given a structure of an affine variety étale over $\Til D$ which makes it the affine envelope of $X$.
\et
A tentative of proof.
We have 
$$
{\Til D}=\k^n\smi\{P=0\}
$$
$$
\ocal (\Til D)=\k[x_1,\ldots,x_n,P,1/P]
$$
with $P\in\k[x_1,\ldots,x_n]$.

We have the diagram
\begin{equation}
  \xymatrix{
  X \ar@{^{(}->}[r]^{\chi}\ar[d]_p &\Smfr(X)\ar[d]^{\Til p}\\
  D\ar@{^{(}->}[r]^{\sf j} &\Til D}
 \end{equation}
where $\Til p$ is the map $h\mapsto h\circ p^\ast\circ {\sf j}^\ast$: $\Til p$ is onto.

Since $P$ is invertible in $\ocal(\Til D)$, for a given $h\in\Smfr(X)$ there is $c_h\in\k$ such that $h\big((p^\ast\circ {\sf j}^\ast)(c_hP)\big)=1$. Fix $h\in\Smfr(X)$ and let $Q_h:=(p^\ast\circ {\sf j}^\ast)(c_hP)$. Moreover, since $p$ is an étale morphism, for any $x\in X$ the homomorphism  
$$
{\Hat p}_{p(x)}:{\Hat\ocal}_{D,\,p(x)}\to\Hat\ocal_{X,x}
$$
is an isomorphism.

Now use Artin theory (\cite{Ar}) to prove first that the functor (\ref{fun1}) is representable by an algebraic space $\Til X$ and then show that $\Til X$ is an affine variety.
\subsection{Other envelopes}\label{oten}
Given an algebraic variety $X$ of pure dimension $n$, let 
\bit
\item[a\rp] ${\mathcal C}^d(X)$ be the set of all connected pure $d$-dimensional subvarieties of $X$
\item[b\rp] $\dcal(X)$ be the set of the Cartier divisors of $X$.
\eit
Given a point $x\in X$, denote ${\mathcal C}^d(X,x)$ the subset of those elements of ${\mathcal C}^d_X$ passing through $x$ and $\dcal(X,x)$ that of those elements of $\dcal(X)$ whose support contains $x$.  

For every subset $K$ of $X$ we consider the set ${\rm hull}^d(K)$ of the points $x$ characterized by the following property: if $
 Z\in {\mathcal C}^d_{X,x}$ and $x\in Z$, then $ Z\cap K\neqv\ES$. The set ${\rm hull}^d(K)$ is called the ${\mathcal C}^d_X$-{\rm envelope} of $K$.
 
 Replacing ${\mathcal C}^d_{X,x}$ by $\dcal(X,x)$ and the above property by ``if $\s\in \dcal(X,x)$ and $x\in {\rm supp}\,\s$, then $ {\rm supp}\,\s\cap K\neqv\ES$'' we get the notion of $\dcal(X,x)$-{\rm envelope} of $K$.  
 
 \vskip1.4truecm
\section{Convexity and concavity}\label{cxcv}
\subsection{Convexity}
The first notion of convexity in algebraic geometry goes back to Goodman and Landman \cite{GL} who define a non-complete algebraic variety $X$ {\em algebraically convex} or (A-{\em convex}) if the following property
\bit
\item[{\rm AC}:]for every open embedding $X\hra Z$ into an algebraic variety $Z$ and for every $x\in{\rm b}_Z X:=Z\ssmi X$ there exists a function $f\in{\mathcal O}(X)\sbs{\mathcal R}(Z)$ which is non-regular at $x$
\eit
is fullfilled.

If $X$ is normal and algebraically convex, then, for every open immersion $X\hra Z$ as above, $Z\ssmi X$ is an algebraic subvariety of codimension one at every point (use the normalization of $Z$).
 
 \bl\label{algco}
 An algebraic variety $X$ is A-convex if and only if the condition {\rm AC} holds for all completions of $X$. If $X$ is normal, condition {\rm AC} for one normal completion suffices.
 \el
$X\hra Z$ and let $Z\hra \oli Z$ be a completion of $Z$. Let $x_0\in {\rm b}_Z X$. Since ${\rm b}_Z X\sbs{\rm b}_{\oli Z}X$ by the hypothesis there exists $f\in{\mathcal O}(X)\sbs{\mathcal R}(\oli Z)$ which is non-regular at $x_0$. This shows that the condition AC holds for any open embedding $X\hra Z$. 

Other possible conditions of algebraic convexity are the following.                                                                                                                                                                                                                                                                                                                                                                                                                                     Let $X\hra Z$ be an arbitrary open immersion. Consider the following conditions:
\bit
\item[\rm AC1:] for every $z_0\in{\rm b}_Z X$ there is $f\in \ocal(X)$ with $z_0\in (f)_\IN$ and $z_0\notin (f)_0$\footnote{$(f)_\IN$ and $(f)_0$ are the poles and the zeros of $f$, respectively.}; equivalently, $f(x_\l)\to\bm{\IN}$ for every net $\{x_\l\}_{\l\in \Lambda}\sbs \oli X\ssmi X$ converging to $z_0$; \item[\rm AC2:] for every $z_0\in{\rm b}_Z X$ and every irreducible curve $C$ through $z_0$ there exists $f\in{\mathcal O}(X)$ such that 
$$
C\nsubseteq (f)_0,\>\>x_0\in (f_{|C})
_\IN.
$$
\eit
$X$ is said {\em strongly $\ocal$-convex} if satisfies the condition AC1 and {\em weakly $\ocal$-convex} if it satisfies the condition AC2.

 If $X$ is strongly $\ocal$-convex then it is A-convex and weakly $\ocal$-convex.
 \br\label{VV}
 The conditions AC1, AC2 can be generalized to coherent sheaves using the following theorem proved by Van Oystaeyen and Verschoren (see \cite[3]{VaVe}.): if $X\hra Z$ is an open immersion of algebraic varieties, the natural morphism
$$
r_X^Z:{\sf Coh}(Z)\longrightarrow {\sf Coh}(X)
$$
is onto.
\er

Each of the conditions {\rm AC1} and {\rm AC2} holds for every open embedding $X\hra Z$ if and only if it holds for one particular completion.
 
 A strongly $\ocal$-convex complex algebraic variety $X$ is a holomorphically convex complex space. The example quoted in Remark \ref{rem} shows that, in general, the converse is not true.

 The following proposition will be used in the sequel
 \bp\label{35p} Let $f:X'\to X$ be a proper morphism of algebraic varieties. Then
\bit
 \item[{\rm 1}\rp] if $X$ is {A}-convex {\lp}strongly $\ocal$-convex, weakly $\ocal$-convex, semi-affine{\rp} then $X'$ is {A}-convex {\lp}strongly $\ocal$-convex, weakly $\ocal$-convex, semi-affine{\rp}. In particular, the normalization $X^{\rm nor}$ of a strongly $\ocal$-convex variety $X$ is strongly $\ocal$-convex;
\item[{\rm 2}\rp]  if $f$ is a modification\footnote{i.e. $Y'\neqv X'$, $Y\neqv X$ are closed subvarieties. $p$ and $f$ are proper, surjective morphisms, $f$ is an isomorphism $X'\ssmi Y'\simeq X\ssmi Y$ and if $x\in Y'$ then $f$ is not a local isomorphism at $x$. $Y'$ is called the {\em exceptional subspace} of the modification, $Y$ the {\em center}. If $f:X'\to X$ is a modification as above with $X$ normal, then, as a consequence of the ZMT, we have the following: ${\rm codim}_\Bbbk{Y_{p(x')}}>1$ at every point $x\in Y'$ and $f^*:\ocal(X)\to\ocal(X') $ is an isomorphism of $\k$-algebras. If $Y$ is a point, $f$ is said a {\em point modification. }}, $X'$ is strongly $\ocal$-convex and $X$ is normal, then $X$ is strongly $\ocal$-convex.
\eit
\ep
 \demo
1) Let $ \oli X$, $\oli {X'}$ be completions of $X$, $X'$ respectively and $x'\in{\rm b}_{\oli {X'}}X'$. Since $X$ is A-convex there exists a net $\sigma=\{x_\l\}_{\l\in\L}$ converging to $x=f(x')$ and $h\in \ocal(X)$ such that $h(x_\l))\to\IN$. Let $\oli{\sigma'}=\oli{f^{-1}(\sigma)}$. If $\oli\sigma'\cap{\rm b}_{\oli X'} X'=\ES$ then $\oli\sigma'$ is complete subvariety of $X'$, therefore, since $f$ is proper, $f(\oli\sigma')=X'\times_X \sigma$ is complete subvariety of $X$, hence $f(\oli\sigma')\cap{\rm b}_{\oli X} X=\ES$: contradiction. Therefore, along a subnet of $\s'$, the function $h\circ f\in \ocal(X')$
is divergent. This shows that $X'$ is A-convex.

The same proof runs as well as for the other two cases, $f$ being proper.  

\vspace{0.5cm}
 As for 2) let
\begin{equation}
  \xymatrix{
    Y' \ar@{^{(}->}[r]^{}\ar[d]_p & X' \ar[d]^f \\
    Y  \ar@{^{(}->}[r]^{} & X 
  }
\end{equation}
be a modification, $X\hra \oli X$ and $X'\hra \oli {X'}$ completions of $X$ and $\oli {X'}$ respectively. Let $\s=\{x_\l\}_{\l\in\L}$ be a net in $X$ converging to a point $x\in {\rm b}_{\oli X} X$. Lift $\s$ by $f$ and let $\s'=\oli{f^{-1}(\s)}$. If $\s'\cap{\rm b}_{\oli X}X \neq\ES$, we may suppose that $\s'$ converges to a point $x\in\oli {X'}\ssmi X'$. Then, by hypothesis, there exists $h\in\ocal(X')$ such that $h\to\bm{\IN}$ along $\s'$. Since $X$ is normal, $h=h_1\circ f$ with $h_1$ regular in $X$. It follows that $h_1\to\bm{\IN}$ along $\s$ and this shows that $X$ is A-convex. If $\s'\cap{\rm b}_{\oli X} \neq\ES$, then $\s'$ is complete and therefore $f(\s')$ is also complete since $f$ is proper. This is absurd since $f(\s')\sps\s$. 
\enddemo 
One of the main results proved by Goodman and Landman in the quoted paper \cite{GL} is the following (see \cite[Corollary 4.4]{GL}): a semi-affine algebraic variety $X$ is A-convex. Conversely, if $X$ is A-convex and $\ocal(X)$ is noetherian, then $X$ is semi-affine. 

This result generalizes in the following way:
\bt\label{p35}
An algebraic variety $X$ is semi-affine if and only if
\bit
\item[\rm i\rp] $\Smfr_0(X)\equiv\Smfr(X)$; 
\item[\rm{ii}\rp] $X$ is strongly $\ocal$-convex. 
\eit
\et
\demo
\vskip-2mm
Let $X$ be semi-affine, $\pi:X\to X_0$ be an affine reduction of $X$ and $X_0\hra\Bbbk^n$ a closed immersion where 
$$
\Bbbk^n=\big\{[z_1,\ldots,z_n,z_0]\in \P^{n+1}:z_0\neq 0\big\}.
$$
Then $\ocal(X_0)\stackrel{\pi^\ast}\simeq \ocal(X)$. Since $\ocal(X_0)\simeq\k[x_1.\ldots,x_n]/I$, where $I$ is an ideal, and $\ocal(X_0)\simeq\ocal(X)$, Hilbert Nullstellensatz implies i).
 
Let $a=[a_1,\ldots,a_n,0]\in {\rm b}_{\oli {X_0}} X_0$. For some $j$ the rational function $\phi_j=(z_j/z_0)_{|\oli X}$ is regular on $X_0$ and $a\in (\phi_j)_\IN$, $a\notin (\phi_j)_0$. This shows that $X_0$ is strongly $\ocal$-convex and, by Proposition \ref{35p}, that $X$ is strongly $\ocal$-convex. 

Conversely, suppose that i) and ii) are fulfilled. Then by i) the topology on $\Smfr_0(X)$ is determined by the function algebra $\Hat{\ocal(X)}$  (cfr. Section \ref{gen}). 

Let us prove that the natural map $\d:X\to\Smfr_0(X)$ is onto. 

We first show that $\d$ is closed. 

Let $C\sbs X$ be a closed subset, $M^0\in\oli{\d(C)}$ and $\sigma=(x_\l)_{\l\in\Lambda}$ a net in $C$ such that $M_{x_\l}=\d(x_\l)\to M^0$. If $x_\l\to x_0\in X$, then $x_0\in C$ and $M^0=\d(x_0)\in\d(C)$. Otherwise we consider a completion $\oli X$ of $X$ which satisfies the property {\rm AC1}. Then there exists a subnet $\sigma'=(x_{\mu})_{\mu\in\Lambda'}$ of $\s$ which converges to a point $y_0\in{\rm b}_{oli X}X$. Let $f\in\rcal(\oli X)$, regular on $X$ and such that $y_0\in (f)_\IN$ and $y_0\notin (f)_0$. We have $\Hat f(M_{x_\mu})=f(x_\mu)$ and (since $y_0$ is a pole and not an indetermination point) eventually
$f(x_\mu)\not=0$. It follows
$$
1=\Hat f(M_{x_\mu})f({x_\mu})^{-1}\lra \Hat f(M^0)\cdot 0=0
$$
contradiction. Therefore $\d$ is a closed map.

Suppose that $M\in \Smfr_0(X)\ssmi\d(X)$. Since $\d(X)$ is closed, there exists an open neighborhood $W$ of $M$ such that $W\cap\d(X)=\ES$. By definition
$$
W=\Hat g_1^{\,-1}(U)\cap\cdots\cap\Hat g_m^{\,-1}(U),
$$
where 
$$
U=\Bbbk\smi F, \>F=\{a_1\}\cup\cdots\cup \{a_r\}
$$
and $g_1,\ldots,g_m$ are non-constant regular functions in $X$.
Let for $1\le j\le m$, $1\le k\le r$
$$
X_{jk}=\big\{x\in X:g_j(x)=a_k\big\}
$$ 
and $Y=\bigcup_{j,k} X_{jk}$. Should $W\cap\d(X_{jk})$ be empty, there is some $g_j$ such that 
$$
\Hat g_j(\d(X\ssmi Y))=g_j(X\ssmi Y)\sbs F.
$$ 
This is a contradiction. 

Thus, $\d(X)=\Smfr_0(X)$. In particular, if the regular functions $f_1,\ldots,f_q\in\ocal(X)$ have no common zero then they generate $\ocal(X)$.
\enddemo
Another consequence of ii) is the following: for any non-complete algebraic curve $C\sbs X$ there exists $f\in\ocal(X)$ such that $f_{|C}$ is non-constant. Indeed let $\oli X$ be a completion of $X$ and $x_0\in\oli C\cap\oli X$, a regular function $f\in\ocal(X)$ such that $x_0\in (f)_\IN$ and $x_0\notin (f)_0$ satisfies the above condition.  

Condition {\rm i)} alone does not ensure that X is semi-affine {\lp}e.g. take $X={\Tilde{\Bbbk^2}}\ssmi\{a\}$ where $\Tilde{{\Bbbk^2}}$ is the blow-up of $\Bbbk^2$ at the origin and $a$ is a point of the exceptional set{\rp}.  
\bc
For an algebraic variety $X$ the following properties are equivalent: 
\bit
\item[1\rp] $X$ is affine;
\item[2\rp] $X$ is strongly $\ocal$-convex, $\Smfr_0(X)\equiv\Smfr(X)$ and the ring $\ocal(X)$ separates points.
\eit
\ec
\bt\label{36t}
An algebraic variety $X$ is affine if and only if the natural map $\d:X\to\Smfr(X)$ is bijective.
\et
Let $X\in{\sf Var}/\k$ be fixed and ${\rm F}_X$ be the contravariant functor 
\be\label{funct}
{\sf Var}/\k\lra\{\k-{\sf algebras}\}, 
\ee
defined by
\be\label{fun1}
Y\rsa {\rm Hom_{alg}}(\ocal(X),\ocal(Y)).
\ee
Theorem\ref{36t} can be rephrased in this way
\bp\label{fun2}
$X$ is an affine variety if and only if, for every algebraic variety $Y$, the natural map 
$$
\Phi_Y:{\rm Mor}(Y,X)\lra{\rm Hom_{\rm alg}}(\ocal(X),\ocal(Y))
$$
is a bijection, i.e. ${\rm F}_X$ is representable in ${\sf Var}/\k$. 
\ep

We have the following extension theorem:

\bt 
Let $X,Z$ be algebraic varieties and $Y\sbs X$ a subvariety of codimension $\ge 2$ at every point. Assume that $X$ is normal and that $Z$ is non-singular or, more generally, locally factorial $K$-complete and strongly $\ocal$-convex. Then the restriction map 
$$
{\rm Mor}\,(X,Z)\to {\rm Mor}\,(X\ssmi Y,Z)
$$
is a bijection.
\et
In the proof of Theorem \ref{p35} the proof of connectness is necessary.

      
\subsubsection{Algebraic Levi problem} 
 Let $X$ be a Stein space, $U\sbs X$ an open subset. A long standing problem in complex analysis, the so called {\em Levi problem}, asks the following: if $U$ is locally Stein is then $U$ a Stein space?

The Levi problem has been completely solved for manifolds of arbitrary dimension and for singular spaces of dimension 2 (\cite{AN}).

The algebraic analogue of the Levi problem, merely obtained replacing ``Steineness'' by ``affineness'' has an immediate positive answer because then the embedding $U\hra X$ is an affine morphism and consequently $U$ is affine (cfr. \cite[p. 128]{Haa}).
%
The {\em algebraic Levi problem} could then be formulated as follows: is a locally A-convex (respectively strongly $\ocal$-convex, weakly-$\ocal$-convex) open subset $U$ of an A-convex (respectively strongly $\ocal$-convex, weakly-$\ocal$-convex,) algebraic variety $X$ an A-convex(respectively strongly $\ocal$-convex, weakly-$\ocal$-convex) algebraic variety?
                                                                                                                                                                                                                                                                                                                                                                                                                                        \subsection{Concavity}
The notion of concavity concerns the algebraic varieties $X$ of pure dimension $\ge $2. 
%

Such a variety $X$ is called {\em algebraically concave} if there exists a completion $X\hra \oli X$ such that ${\rm b}_{\oli X}X$ is of codimension $\ge 2$. 

$X$ is called {\em strongly $\ocal$-concave} if there is a completion $X\hra \oli X$ such that ${\rm b}_{\oli X}X$ is the excepional subspace of a point modificaion of $\oli X$. 
An affine variety is neither strongly $\ocal$-concave nor algebrically concave. The complement of a finite subset of a complete variety of pure dimension $\ge 2$ is strongly $\ocal$-concave.

As a trivial consequence of the normalization we have
\bp\label{p25}
If $X$ is a $\ocal$-concave variety then $\ocal(X)=\k$.
\ep
%
%
\bp\label{const}
Let $X\hra \oli X$, $X\hra Z$ be two normal completions of the algebraic variety $X$. If ${\rm b}_{\oli X}X$ has a strongly $\ocal$-convex {\rm(}semi-affine{\rm)}neighborhood then the same holds for ${\rm b}_ZX$.  
\ep
\demo
Let $U$ be a strongly $\ocal$-convex neighborhood of ${\rm b}_{\oli X}X$. 
The diagram of maps
\begin{equation}
\begin{split}
\xymatrix{X\>\ar@{^{(}->}[r]\ar@{^{(}->}[d]& \oli X\\
Z& }
\end{split}
\end{equation}
defines a birational morphism $\r:\oli X\to Z$. If $\Gamma_\r$ is the graph of $\r$, then 
\be\label{mod2}
 \xymatrix{
p_{\oli X}^{-1}({\rm b}_{\oli X}X)\cap\Gamma_\r\ar@{^{(}->}[r]\ar[d]_{p_{\oli X}}&\Gamma_\r\cap p_{\oli X}^{-1}(U)\ar[d]_{p_{\oli X}} \\
    {\rm b}_{\oli X}X\ar@{^{(}->}[r] &U 
  }
  \ee
 \be\label{mod3}
 \xymatrix{
p_{\oli X}^{-1}({\rm b}_{\oli X}X)\cap\Gamma_\r\ar@{^{(}->}[r]\ar[d]_{p_{\oli Z}}&\Gamma_\r\cap p_{\oli X}^{-1}(U)\ar[d]_{p_{\oli Z}} \\
    {\rm b}_{\oli Z}X\ar@{^{(}->}[r] &\Sigma_\rho(U) 
  }
  \ee
are modifications. Since $U$ is normal and $\ocal$-convex, Proposition \ref{35p} applies: $\Gamma_\r\cap p_{\oli X}^{-1}(U)$ is strongly $\ocal$-convex, therefore, by part 2) of the same Proposition, $\Sigma_\rho(U)$ (which is a neighborhood of $X$ in $Z$) is strongly $\ocal$-convex. 
\enddemo
\bt\label{t2}
Let $X$ be a normal strongly $\ocal$-concave variety of pure dimension $\ge 2$, $\mathcal F\in {\sf Coh}(X)$ be a torsion free, coherent sheaf on $X$. Under the condition of the last proposition $\Gamma (X,\mathcal F)$ is a finite dimensional $\k$-vector space. 
\et
\demo
Let $X\hra \oli X$ be a normal completion of $X$ and $U$ a special neighborhood of ${\rm b}_{\oli X}X$.
$\fcal$ extends to $\oli X$ by a torsion free coherent sheaf $\oli\fcal$. Indeed, by \cite[3.1 Proposition]{VaVe}, there exists $\Hat\fcal\in {\sf Coh}(\oli X)$ which extends $\fcal$ to $\oli X$. If $\tcal(\Hat\fcal)$ denotes the torsion of $\Hat\fcal$ the $\oli\fcal:=\Hat\fcal/\tcal(\Hat\fcal)$ is the desired extension. 

Let $U$ be a neighborhood of ${\rm b}_{\oli X}X$ and
 $f:U\to V$ be a point modification with $V$ affine  and $V$ normal: $f_*\oli\fcal$ is also torsion free. Let $s\in\Gamma(X)$ be a global section of $\fcal$. Then $\s:=s_{\vert U}$ defines a section $\s_*\in\Gamma(V\ssmi f({\rm b}_{\oli X}X),\ocal_V)$ which extends uniquely to a section $\oli\s\in\Gamma(V,f_*\oli\fcal)$ since $V$ is normal, the sheaf $f_*\oli\fcal$ is torsion free and $f({\rm b}_{\oli X}X)$ is finite. It follows
\[
\kdim \Gamma(X,\fcal)=\kdim \Gamma(\oli X,\oli\fcal)<+\infty
\]
and this concludes the proof.
\enddemo
\subsubsection{$d$-concavity.}
A more geometric definition of concavity is based on the notion of ``envelope''. Let ${\ccal}^d_X$ be the set of all connected,
complete and $d$-dimensional algebraic subvarieties of an algebraic variety $X$. Given a point $x\in X$, denote ${\ccal}^d_{X,x}$ the
subset of those elements of ${\ccal}^d_X$ passing through $x$. Let $K$ be a subset of
$X$. We define the ${\ccal}^d_{X}\!-\!envelope$ of $K$ in $X$ by 
$$
{\rm Env}l_d(K)=\{x\in X:Z\cap K\neqv\ES \>\fa Z\in {\ccal}^d_{X,x} \}.
$$
An algebraic variety $X$ is said to be $d$-{\em concave} if there exists a
complete $(n-d)$-dimensional subvariety $Y$ such that ${\rm  Env}_d(Y)$ has a nonempty
interior part.

A complete variety is $d$-concave for every $d\ge 0$. An affine variety of positive dimension is never $d$-concave.

Here are some natural questions  
\begin{enumerate}
\item[1\rp] Is the complement $\C{\P}^{\,n}\ssmi Y_d$, of a subvariety $Y_d$ of pure 
dimension $d$, $(d+1)$-concave ?\\
\item[2\rp] Let $X$ be a complete algebraic variety of dimension $n$ and $Y_d$ a subvariety
of pure dimension $d$. Assume that $X\ssmi Y_d$ is $q$-concave with $q<d+1$. Is then
$X$ contractible along $Y_d$ in the category of the algebraic spaces?
\item[3\rp] Let $X$ be a $1$-concave algebraic variety and $\oli X$ one completion of $X$.
Does a neighborhood of $\oli X\ssmi X$ have some convexity?\\
\end{enumerate}
\section{Convex and concave line bundles}
\subsection{Convex line bundles}
\nin Let ${\sf L}\stackrel{\pi}{\rightarrow}X$ be a line bundle on the algebraic variety $X$, $\ocal(\sf L)$ the sheaf of the germs of regular sections of $\sf L$ and ${\sf L}^\ast$ the dual of $\sf L$. We denote $\oli{\sf L}\stackrel{\oli\pi}{\rightarrow}X$ the projective closure $\P({\sf L}\oplus\k)$ of ${\sf L}\stackrel{\pi}{\rightarrow}X$.
We say that ${\sf L}\stackrel{\pi}{\rightarrow}X$ is 
 {\em strongly $\ocal$-convex along the fibres} if for every point $z\in\oli{\sf L}\ssmi{\sf L}$ and a net $\sigma=(z_\l)_{\l\in I}$ on ${\sf L}$ converging to $z$ there exists a regular function $f:{\sf L}\to\Bbbk$ which is divergent on $\sigma$. Observe that, if $X$ is complete and ${\sf L}$ is strongly $\ocal$-convex along the fibres then (the total space of the line bundle) ${\sf L}\stackrel{\pi}{\rightarrow}X$ is strongly $\ocal$-convex variety.
 
More generally 
 \bt\label{t9}
Let $X$ be strongly $\ocal$-convex variety and ${\sf L}\stackrel{\pi}{\rightarrow}X$ a line bundle which is strongly $\ocal$-convex along the fibres. Then, ${\sf L}$ is a strongly $\ocal$-convex variety. 
\et
 Now assume that $X$ is complete and $\{{\sf L}_{\vert U_i}\}$ is a trivialization of $\sf L$ with respect to an affine covering $U_1, U_2n\cdots,U_m$ of $X$. Let
$\{h_{ij}\}_{ij}$ denote the corresponding cocycle and $\z_j$ the fiber coordinates. 

If $f$ is a regular function in $\sf L$ then 
$$
f_{\vert \pi^{-1}(U_j)}=\sum_{r=1}^{m_j}f_{r_,j}\z^r_j, f_{r_,j}\in\ocal(U_j)
$$ 
where the functions $f_{r_,j}$ are uniquely determined by $f$. On $\pi^{-1}(U_j)\cap \pi^{-1}(U_l)(\neqv\ES)$ it results $f_{0,j}=f_{0,l}$, $m_j=n_l$ and the functions $f_{r_,j}$ define a regular section $\s_r$ of the dual ${\sf L}^{\ast r}$ of ${\sf L}^{\ast r}$. Conversely, if $\s_r=\{f_{r_,j}\}_j$ is a section of ${\sf L}^{\ast r}$, $f_{r,j}\z^r_j$ is a regular function on $L$. It follows that 
$$
\s_0\oplus \s_1\cdots\oplus \s_m\mapsto f_{0_,j}+f_{1,j}\z_j\cdots+f_{m,j}\z^m_j
$$
define a surjective homomorphism of $\k$-algebras
\be\label{gr}
\acal(X,{\sf L}^\ast)\lra\G(\sf L,\ocal_{\sf L})
\ee
where
$$
 \acal(X,{\sf L}^\ast)=\bigoplus_{r\ge 0}\acal_r(X,{\sf L}^\ast):=\bigoplus_{r\ge 0}\G(X,{\sf L}^{\ast r}).
$$
We want to state a criterion for finite generation of $\acal(X,{\sf L}^\ast)$. It is based on the notion of a-sequence.
A system $s_1,\ldots, s_m$ of sections of ${\sf L}^{\ast r}$, $r\ge 1$, is called an {\em-a-sequence of length $r$} for $\acal(X,{\sf L}^\ast)$ if there exists an integer $h_0$ such that
\beqn
&&\big(\acal_h(X,{\sf L}^\ast) s_1+\ldots+\acal_h(X,{\sf L}^\ast) s_i\big)\cap \acal_h(X,{\sf L}^\ast)s_{i+1}=\\
&&\big(\acal_{h-r}(X,{\sf L}^\ast) s_1+\ldots+\acal_{h-r}(X,{\sf L}^\ast) s_i\big)s_{i+1}
\eeqn
for all $i$ with $1\le i\le m-1$ and for $h\ge h_0$.
\bp\label{fin}
Let ${\sf L}$ be a line bundle on a complete $n$-dimensional algebraic variety $X$. If $\G(X,{\sf L}^{\ast r})$, for some $r\ge 1$, contains $n+1$ sections that give an a-sequence for $\acal(X,{\sf L}^\ast)$, then the function algebra $\acal(X,{\sf L}^\ast)$ is finitely generated.  
\ep
\demo
We know that 
$$
\kdim \G(X{\sf L}^{\ast r})\le Cr^h, h\le n, \forall r   
$$
 for some positive $C$ (see  \cite{A}), therefore Proposition 7 of \cite{AT} applies.
\enddemo
Sufficient conditions for the existence of a-sequences are established in \cite{AT}.
\bt\label{t6}
Under the hypothesis of Proposition {\rm\ref{fin}}, assume that ${\sf L}$ satisfies the following additional condition: for every $x\in X$ there exists a section $s\in\Gamma(X,{\sf L}^*)$ such that $s(x)\neq 0$. Then
\bit 
\item[\rm1)] ${\sf L}$ is a semi-affine variety, $f:{\sf L}\to{\sf L}_0$;
\item[\rm2\rp] if $\kdim{\sf L}=\kdim{\sf L_0}$, ${\sf L}$ is a point modification of an affine variety $Z$ with the 0-section of ${\sf L}$ as exceptional set. 
\eit
\et
\demo
1) Since $(X,{\sf L})$ is finitely generated, ${\sf L_0}={\sf L}$, therefore, taking into account Theorem \ref{p35}, we have to prove that ${\sf L}$ is strongly $\ocal$-convex.
Let $\sigma=\{z_\l\}_{\l\in I}$ be a net on ${\sf L}$ converging to $z\in\oli{\sf L}\ssmi{\sf L}$ and $x=\oli\pi(z)$. Let $s\in\Gamma(X,{\sf L}^*)$ be such that $s(x)\neq 0$, $s=\{s_i\}_i$ with respect to a trivialization of ${\sf L}^*$. Then, if $\zeta_i$ denotes the generic fiber coordinate in the given trivialization, the regular function $h=\{s_i\zeta_i\}$ diverges on $\sigma$ ($f(x_\l)\to\bm{\IN}$). This shows that ${\sf L}$ is strongly $\ocal$-convex along the fibres, therefore it is strongly $\ocal$-convex, $X$ being complete. This proves part 1).

As for part 2), let $f:{\sf L}\to{\sf L}_0$ be the affine reduction of $\sf L$ and assume that $\kdim{\sf L}=\kdim{\sf L_0}$. We have the following:
\bit
\item[\rm a)]$\sf L$ contains but finitely many of complete subvarieties of maximal dimension.
\eit
Let $S\sbs\sf L$ be a complete subvariety of maximal dimension. Then, since $S\cap\pi^{-1}(x)$ is finite, for every $x\in X$, $\pi(S)=X$. If $y\in{\sf L_0}$ denotes the point $f(S)$, then, for every $x\in X$, the non-complete curve $C_x=(\pi^{-1}(x)$ contains $y$. If, by a contradiction, there were infinite subvarieties $S$, the points $y$ would be infinite and this is absurd.This proves a).

By the Stein factorization (see \cite[4.3.1]{GD}) we have a commutative diagram of morphisms of algebraic varieties 
$$
\xymatrix{{\sf L}\ar[rr]^{f'}\ar[rd]_{f}& &Z\ar[ld]^{g}\\ &{\sf L_0} &}
$$
where $f'$ is proper, $f'_\star\ocal_{\sf L}=\ocal_{Z}$ and $g$ is finite. By ZMT the fibres of $f'^ {-1}(z)$ are connected. Since ${\sf L_0}$ is affine and $g$ is finite, $Z$ is affine. Moreover, by a), away from a finite subset of $Z$, the map $f'$ is bijective and therefore an isomorphism since $f'_\star\ocal_{\sf L}=\ocal_{Z}$.   
\enddemo
Arguing as in \cite{Gr}, as a corollary we get the following existence theorem for modifications (see also \cite[8.9]{gd}):
 \bt\label{t4}
 Let $Z$ be a non-singular algebraic variety, $X\sbs Z$ a complete non-singular hypersurface. Assume that the normal bundle ${\sf N}_{X/Z}$ satisfies the conditions of Theorem {\rm \ref{t6}} for ${\sf L}$. Then $X$ is exceptional in $Z$. 
 \et
 \subsection{Concave line bundles} 
 Let $s=\{s_i\}_i$ be a regular section of the line bundle ${\sf L}$ on the algebraic variety $X$. If $\xi_i$ is the fiber coordinate, $f^\ast_s:=\{s_i/\xi_i\}_i$ is a rational function in ${\sf L}$, regular away from the zero-section ${\bf 0}_L$. We call the bundle ${\sf L}$  {\em concave along the fibers} if for every $x\in X$ there exists a section $s\in\Gamma(X,{\sf L})$ such that $s(x)\neq 0$.
 
Let us end this section stating without proof the following two theorems.
 
\bt\label{t4}
Let $X$ be a complete, normal algebraic variety. Let ${\sf L}\to X$ be a line bundle concave along the fibres and satisfying the hypothesis of Proposition \ref{fin}. Then ${\sf L}$ is a strongly $\ocal$-concave variety.
\et
\bt\label{t5}
 Let $\fcal$ be a coherent sheaf on a normal concave variety $X$ and $F\to X$ a concave line bundle. Then 
$$
\Gamma(X,\fcal\tens\ocal(L^{*n}))=0. 
$$
for $n\gg0$.
\et


\begin{thebibliography}{HST17b}

\bibitem[A]{A} A. Andreotti, Théorèmes de dépendence algébrique sur les espaces pseudoconcaves, {\em Bull. Soc.  Math. France} \textbf{I91}, (1963), 1--38.

\bibitem[AnNa]{AN} A. Andreotti, R. Narasimhan, Oka's Heftungslemma and Levi Problem for complex spaces, {\em Trans. Amer. Math. Soc.} \textbf{III}, (1964), 256--277.

\bibitem[AnTo]{AT} A. Andreotti, G. Tomassini, Essays on Topology and related Topics. Memoires dédiés à Georges de Rham, 85-104. Publiés sous la direction de André Haefliger et Raghavan Narasimhan, Springer-Verlag 1970. 

 \bibitem[Ar]{Ar} M. Artin,  Algebraic spaces, James K.Whittemore Lectures in Mathematics given at Yale University, 1969. Yale Mathematical Monographs, 3. Yale University Press, New Haven, Conn.-London, 1971. vii+39 pp.
 
 \bibitem[De]{De} J.-P.  Demailly,  Complex Analytic abd Differential Geometry. Université de Grenoble, Institut Fourier, UMR 5582 du CNRS38402 Saint-Martin d'Hères, France. Version of Thursday June 21,  2012.
 
%
\bibitem[GoLa]{GL} J. E. Goodman, A. Landman, Varieties Proper over Affine Schemes, {\em Inven. math.} \textbf {20}, (1973), 267-312.
 
  \bibitem[Gra]{Gr} H. Grauert, \"Uber Modifikationen und exzeptionelle analytische Mengen, {\em Math. Ann.} \textbf {146}, (1967), 331--368.

%
\bibitem[EGA II]{gd} A. Grothendieck, J. Dieudonné, \'Elements de  g\'eom\'etrie alg\'ebrique II {\em Publ. Math. IHES}, \textbf{11} (1961).

\bibitem[EGA III]{GD} A. Grothendieck, J. Dieudonné, \'Elements de g\'eom\'etrie alg\'ebrique III {\em Publ. Math. IHES}, \textbf{11} (1961).

\bibitem[Ha1]{Ha} R. Hartshorne, Ample Subvarieties, {\em Lec. Notes in Math.} \textbf{156}, Springer-Verlag 1970.

\bibitem[Ha2]{Haa} R. Hartshorne, Algebraic Geometry, {\em Graduate Texts in Math.} \textbf{52}, New York: Springer-Verlag 1970.


\bibitem[Mi]{Mi} J. S. Milne, Elliptic Curves, World Science 2020.


\bibitem[Na]{Na} A. Nagata, Imbedding of an abstract variety in a complete variety, {\em J. Math. Kyoto Univ. } \textbf{2}, (1962), 1--10. 

\bibitem[Ne]{Ne} A. Neeman, Affine and Hilbert's Fourteenth Problem, {\em Ann. of Math. } \textbf{127}, (1962), 229--244.

\bibitem[Se]{Se} J. P. Serre, G\'eom\'etrie alg\'ebrique et g\'eom\'etrie analytique, {\em Ann.Inst. Fourier} \textbf{6}, (1955-1956), 1--42. 

\bibitem[Si]{Si}Y. T. Siu, Every Stein subvariety admits a Stein neighborhood, {\em Inv. Math.} \textbf{38}, (1976-1977), 89--100.

%

\bibitem[VaVe]{VaVe} F. Van Oystaeyen, A. Verschoren, Extending Coherent and Quasi-Coherent Sheaves on Generally Closed Spaces, {\em J. Algebra} \textbf{89}, (1984), 224--236. 

\end{thebibliography}
\end{document}